\documentclass[12pt]{amsart}
\usepackage{amssymb}
\usepackage{colordvi}
\usepackage{epsf, epsfig}

\copyrightinfo{}{}
\title{Tropical Interpolation}
\author{Frank Sottile}
\address{Department of Mathematics\\
         Texas A\&M University\\
         College Station\\
         TX \ 77843\\
         USA}
\email{sottile@math.tamu.edu}
\urladdr{http://www.math.tamu.edu/\~{}sottile}

\thanks{Sottile supported in part by NSF CAREER
  grant DMS-0134860 and the Clay Mathematical Institute}
\thanks{We gratefully thank our editor, Silvio Levy and the MSRI
  members whose work we describe.}

\begin{document}

\maketitle

Everyone knows that two points determine a line, and many people who
have studied geometry know that five points on the plane determine a
conic.  In general, if you have $m$ random points in the plane and you
want to pass a rational curve of degree $d$ through all of them,
there may be no solution to this interpolation problem (if $m$ is
too big), or an infinite number of solutions (if $m$ is too small),
or a finite number of solutions (if $m$ is just right).  It turns
out that ``$m$ just right'' means $m=3d{-}1$ ($m=2$ for lines and
$m=5$ for conics). 

A harder question is, if $m=3d{-}1$, {\it how many\/} rational curves of degree
$d$ interpolate the points?  Let's call this number $N_d$, so that 
$N_1=1$ and $N_2=1$ because the line and conic of the previous
paragraph are unique.  It has long been known that $N_3=12$, and
in 1873 Zeuthen~\cite{Ze1873} showed that $N_4=620$.
That was where matters stood until 1989, when Ran~\cite{R89} gave a recursion for these
numbers. 
About ten years ago, Kontsevich and Manin~\cite{KM} used associativity in quantum
cohomology of $\mathbb{P}^2$ to give the elegant recursion
\[
  N_d\ =\   \sum_{a+b=d} N_a N_b \left( a^2b^2\binom{3d-4}{3a-2} - 
           a^3b\binom{3d-4}{3a-1}\right)\ .
\]

The research themes in the MSRI Winter 2004 semester on Topological
Aspects of Real Algebraic Geometry
included enumerative real algebraic geometry, tropical geometry,
real plane curves, and applications of real algebraic geometry.
All are woven together in the unfolding story of this interpolation
problem, a prototypical problem of \textit{enumerative 
geometry}, which is the art of counting geometric figures determined by
given incidence conditions.  Here is another problem: how many lines
in space meet four given lines?
To answer this, note that three lines lie on a unique doubly-ruled 
hyperboloid.
\[
  \begin{picture}(280,117)
   \put(0,0){\includegraphics[height=1.6in]{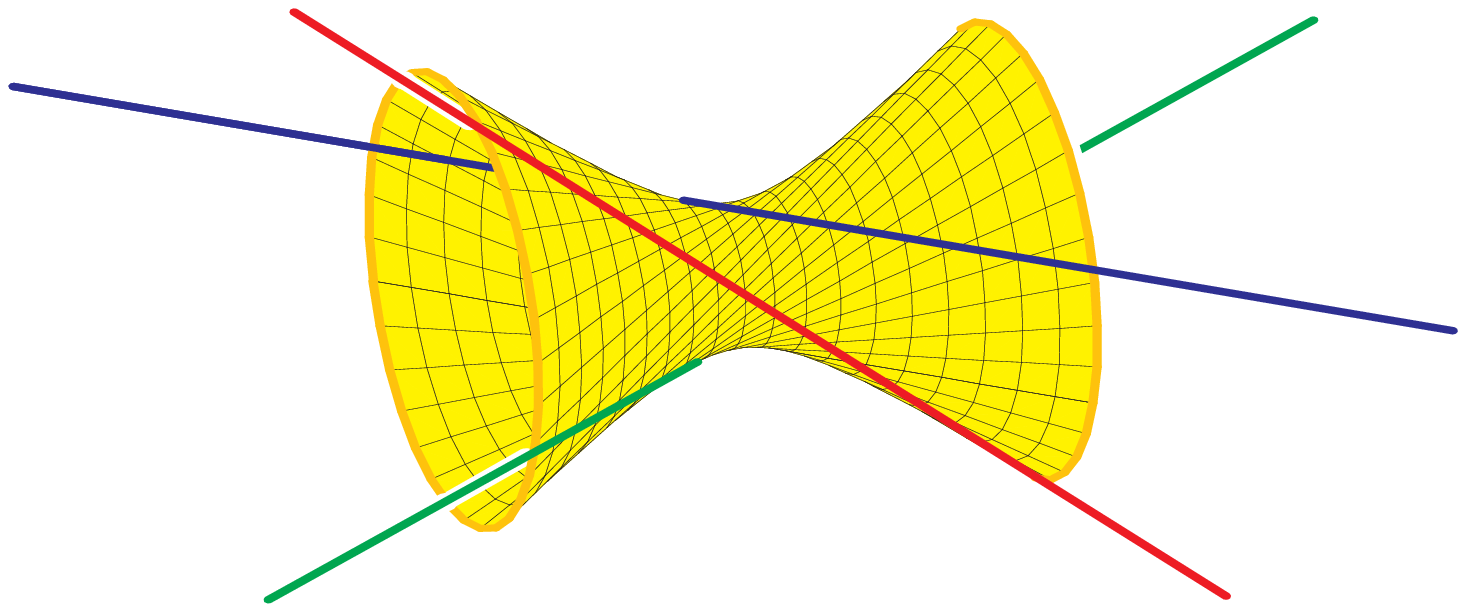}}
  \end{picture}
\]
The three lines lie in one ruling, and the second ruling consists of
the lines meeting the given three lines.
Since the hyperboloid is defined by a quadratic equation, a fourth
line will meet it in two points.
Through each of these two points there is a line in the second
ruling, and these are the two lines meeting our four given lines.

Enumerative geometry works best over the {\it complex} numbers, as
the number of real figures depends rather subtly on the configuration of
the figures giving the incidence conditions.
For example, the fourth line may meet the hyperboloid in two real
points, or in two complex conjugate points, and so there are either
two or no real lines meeting all four.
Based on many examples, we have come to expect that any enumerative
problem may have all of its solutions be real~\cite{So}. 

Another such problem is the 12 rational curves interpolating 8
points in the plane.
Most mathematicians are familiar with the nodal (rational)
cubic shown on the left below.
There is another type of real rational cubic, shown on the right.
\[
  \begin{picture}(250,85)
   \put(0,0){\includegraphics[height=80pt]{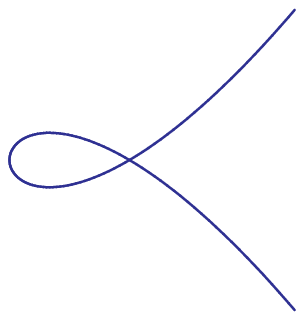}\qquad\qquad\qquad
  \includegraphics[height=80pt]{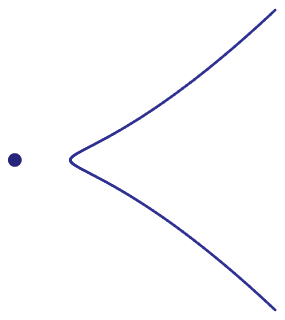}}
 \end{picture}
\]
In the second curve, two complex conjugate branches meet
at the isolated point.
If we let $N(\tau)$ be the number of real curves of type $\tau$ 
interpolating 8 given points, then 
Kharlamov and Degtyarev~\cite{DeKh00} showed that 
\[
   N(\raisebox{-2pt}{\includegraphics[height=10pt]{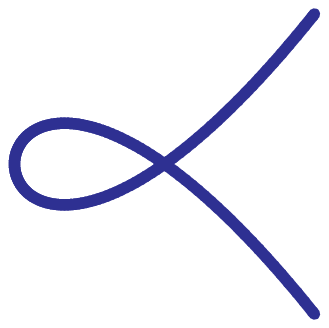}})
   \ -\ 
   N(\raisebox{-2pt}{\includegraphics[height=10pt]{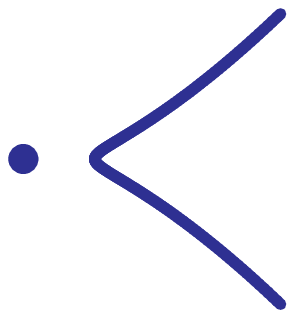}})
    \ =\ 8\,.
\]
Their elementary topological methods are described on the
web$^\dagger$.{\let\thefootnote\relax\footnotetext{$^\dagger$   
{\tt http://www.math.tamu.edu/\~{}sottile/stories/real\_cubics.html}}}

Since there are at most 12 such curves, 
$N(\raisebox{-2pt}{\includegraphics[height=10pt]{figures/Resub.eps}})
 +N(\raisebox{-2pt}{\includegraphics[height=10pt]{figures/Cosub.eps}})
 \leq 12$, 
and so there are 8, 10, or 12 real rational cubics interpolating 8
real points in the plane, depending upon the number 
($0$, $1$, or $2$) of cubics with an isolated point.
Thus there will be 12 real rational cubics interpolating any 8 of
the 9 points of intersection of the two cubics below.
\[
  \begin{picture}(150,75)
   \put(0,0){ \includegraphics[height=70pt]{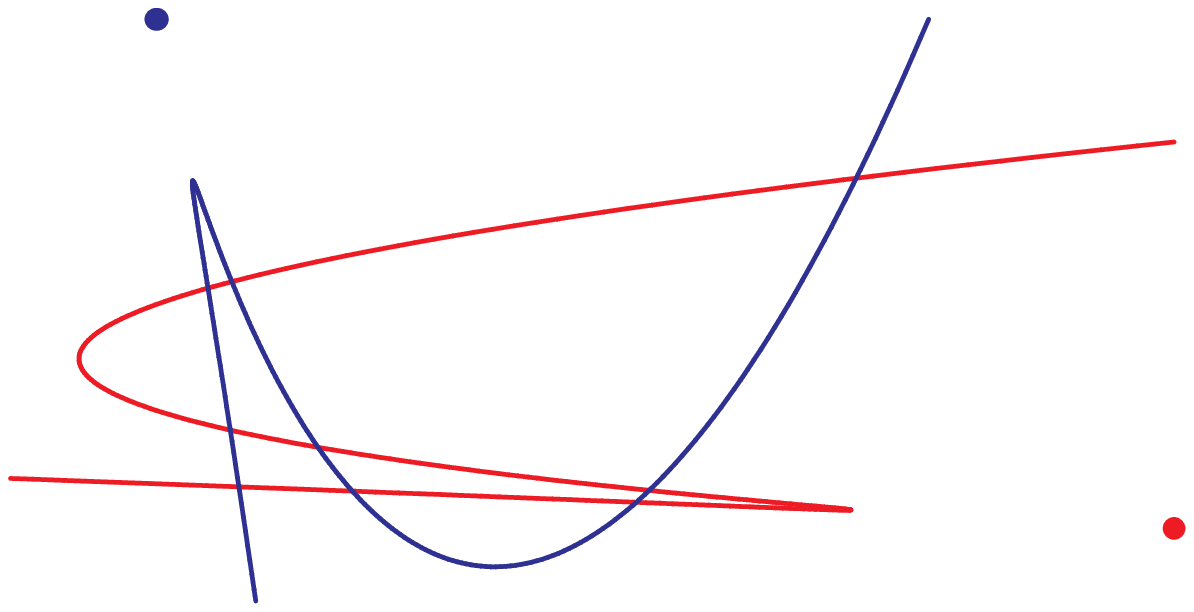}}
 \end{picture}
\]

Welschinger~\cite{W}, who was an MSRI postdoc last Winter, developed 
this example into a theory.
In general, the singularities of a real rational plane curve $C$ are nodes
or isolated points. 
The parity of the number of nodes is its {\it sign}
$\sigma(C)\in\{\pm1\}$.
Given $3d{-}1$ real points in the plane,
Welschinger considered the quantity
\[
  \left|\sum \sigma(C)\right|\,,
\]
the sum over all real rational curves $C$ of degree $d$ that interpolate
the points.
He showed that this weighted sum does not depend upon the choice
of points.
Write $W_d$ for this invariant of Welschinger.
For example, we just saw that $W_3=8$.

This was a breakthrough, as $W_d$ was
(almost) the first truly non-trivial invariant in enumerative real
algebraic geometry.
Note that $W_d$ is a lower bound for the number of real rational curves
through $3d{-}1$ real points in the plane, and $W_d\leq N_d$.
\medskip

Mikhalkin, who was an organizer of the semester, provided the key to
computing $W_d$ using tropical algebraic
geometry~\cite{Mi}. 
This is the geometry of the tropical semiring,
where the operations of $\max$ and $+$ on real numbers 
replace the usual operations of $+$ and $\cdot$.
A tropical polynomial is a piecewise linear function of the form
\[
   T(x,y)\ =\ \max_{(i,j)\in\Delta}\{ x\cdot i + y\cdot j
   +c_{i,j}\}\,,
\]
where $\Delta\subset\mathbb{Z}^2$ is the finite set of exponents of $T$
and $c_{i,j}\in\mathbb{R}$ are its coefficients.
A tropical polynomial $T$ defines a tropical curve, which is the set
of points $(x,y)$ where $T(x,y)$ is not differentiable.
Here are some tropical curves.
\[
  \begin{matrix}
   \includegraphics[width=55pt]{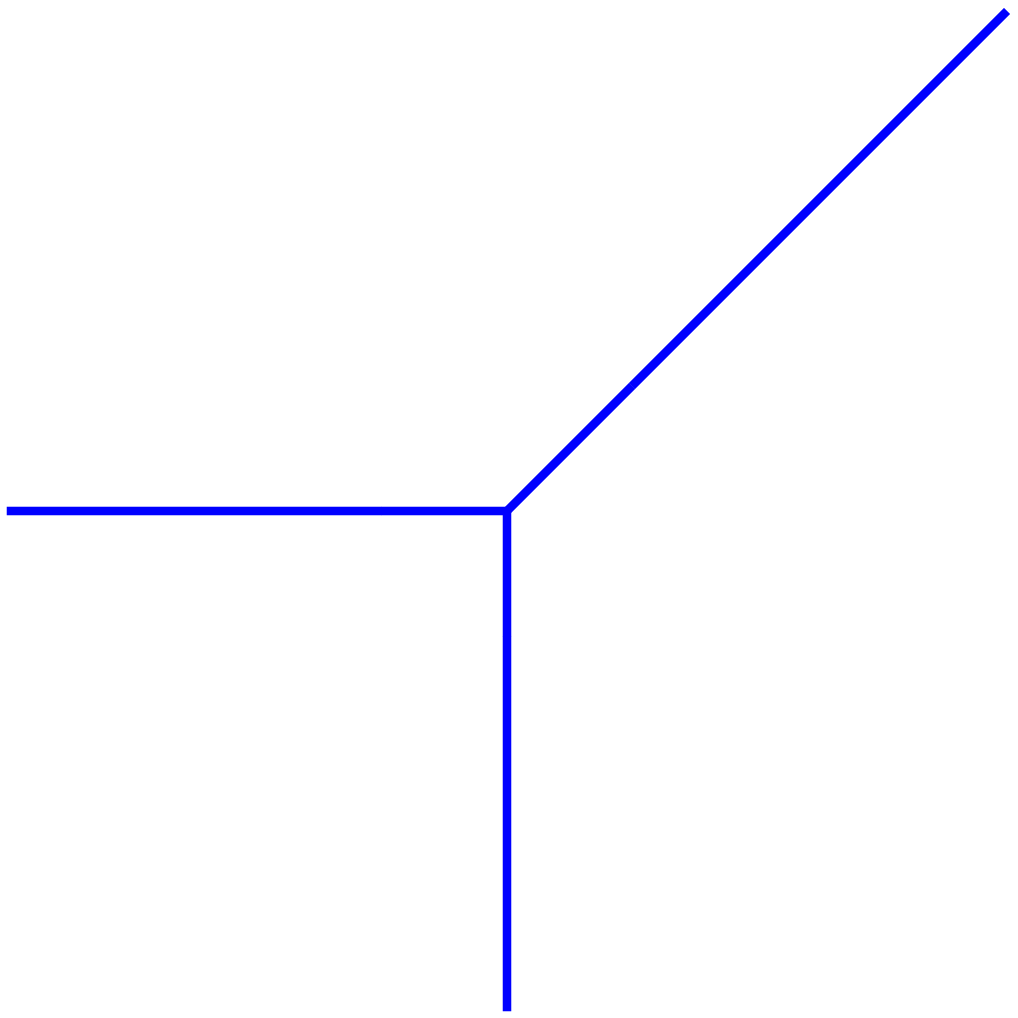}&\quad&
   \includegraphics[width=65pt]{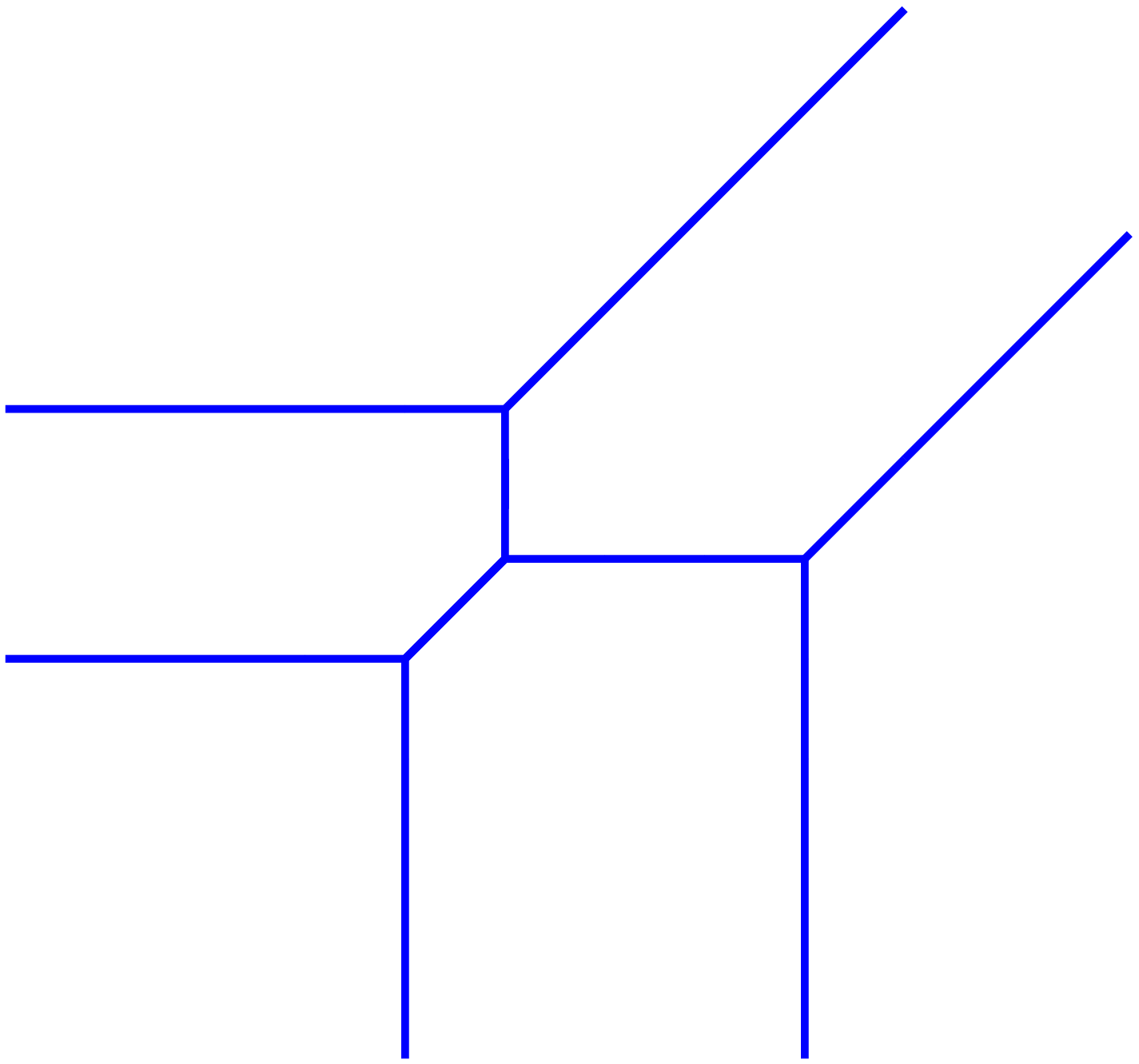}&\quad&
   \includegraphics[width=65pt]{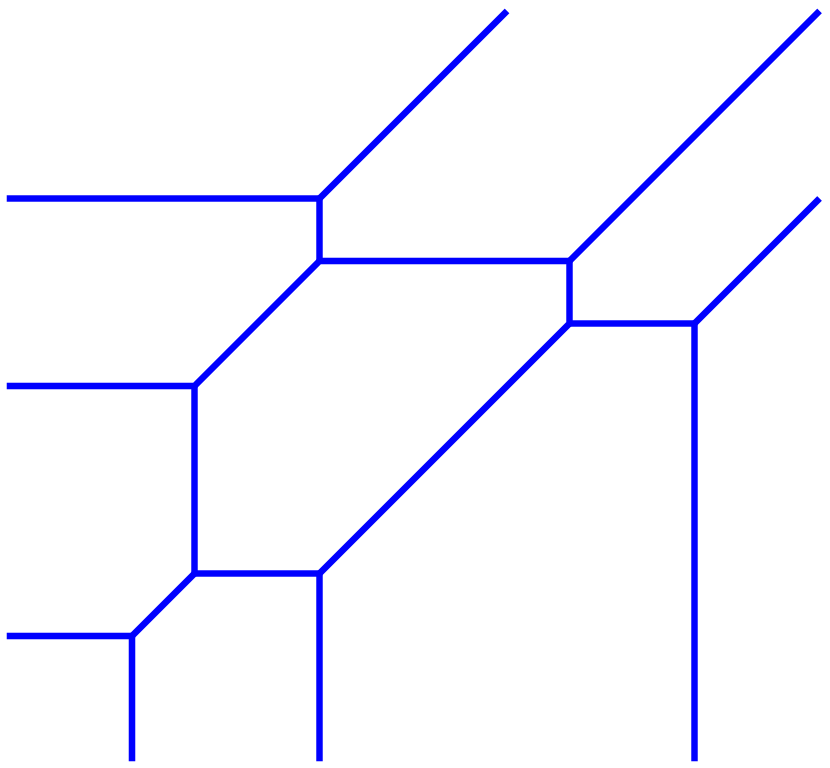}&\quad & 
   \begin{picture}(65,65)
      \put(0,0){\includegraphics[width=65pt]{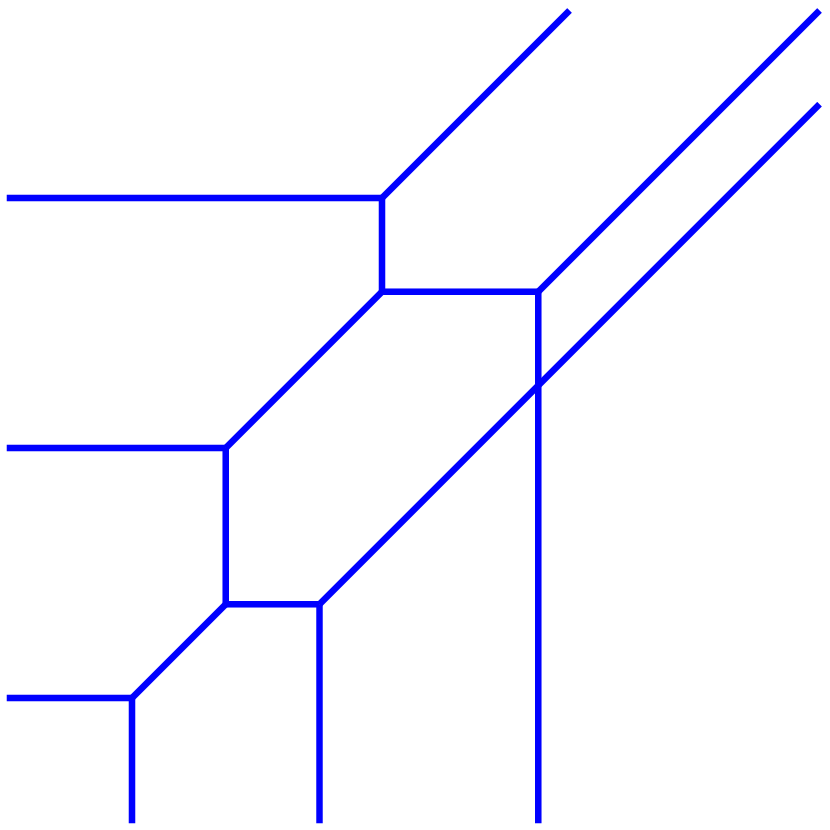}}
      \put(54,23.5){\vector(-1,1){10}}   
      \put(50,16){node}
    \end{picture}\\
  \text{line}&&\text{conic}&&\text{cubic}&&\text{rational cubic}
  \end{matrix}
\]
The degree of a tropical curve is the number of rays tending
to infinity in either of the three directions West, South, or North East.
A tropical curve is rational if it is a piecewise-linear immersion of
a tree. 
Nodes have valence 4.

Mikhalkin showed that there are only finitely many rational tropical curves
of degree $d$ interpolating $3d{-}1$ generic points.
While the number of such curves does depend upon the choice of
points, Mikhalkin attached positive multiplicities to each
tropical curve so that the weighted sum does not, and is in fact
equal to $N_d$.
He also reduced these multiplicities and the enumeration of
tropical curves to the combinatorics of lattice
paths within a triangle of side length $d$.

Mikhalkin used a correspondence involving the map
$\mbox{Log}\colon(\mathbb{C}^*)^2\to\mathbb{R}^2$  
defined  by $(x,y)\mapsto(\log|x|,\log|y|)$, and a certain `large
complex limit' of the complex structure on $(\mathbb{C}^*)^2$.
Under this large complex limit, rational curves of degree $d$ 
interpolating $3d{-}1$ points in $(\mathbb{C}^*)^2$ deform to `complex
tropical curves', whose images under $\mbox{Log}$ are ordinary
tropical curves interpolating the images of the points.
The  multiplicity of a tropical curve $T$ is the number of
complex tropical curves which project to $T$.

What about real curves?
Following this correspondence, 
Mikhalkin attached a real multiplicity to each tropical curve and
showed that if the tropical curves interpolating a given $3d{-}1$
points have total real multiplicity $N$, then there are $3d{-}1$ real
points which are interpolated by $N$ real rational curves of degree $d$.
This real multiplicity is again expressed in terms of
lattice paths.

What about Welschinger's invariant?
In the same way,  Mikhalkin attached a signed weight to each
tropical curve (a tropical version of Welschinger's sign) and 
showed that the corresponding weighted sum equals Welschinger's
invariant.
As before, this tropical signed weight may be expressed in terms of
lattice paths.

During the semester at MSRI, Itenberg,
Kharlamov, and Shustin~\cite{IKS} used Mikhalkin's results to 
estimate Welschinger's invariant.
They showed that $W_d\geq\frac{1}{3}d!$, and also 
\[
   \log W_d\ =\ \log N_d + O(d),\qquad
    \mbox{and}\qquad
    \log N_d \ =\ 3d \log d + O(d)\,.
\]
Thus at least logarithmically, {\it most} rational curves of
degree $d$ interpolating $3d{-}1$ real points in the plane are real.

There are two other instances of this phenomenon of lower bounds,
the first of which predates Welschinger's work. 
Suppose that $d$ is even and let $W(s)$ be a real polynomial of
degree $k(d-k+1)$.
Then Eremenko and Gabrielov~\cite{EG} showed that there exist real
polynomials $f_1(s),\dotsc,f_k(s)$ of degree $d$ whose Wronski
determinant is $W(s)$. 
In fact, they proved a lower bound on the number of $k$-tuples of
polynomials, up to an equivalence.
Similarly, while at MSRI, Soprunova and I~\cite{SS} studied sparse
polynomial systems associated to posets, showing that the number of
real solutions is bounded below by the sign-imbalance of the poset.
Such lower bounds to enumerative problems, which imply the existence of
real solutions, are important for applications.

For example, this story was recounted over beer one evening
at the MSRI Workshop on Geometric Modeling and Real Algebraic
Geometry in April 2004.
A participant, Schicho, realized that the result $W_3=8$ for cubics
explained why a method he had developed always seemed to work.
This was an algorithm to compute an approximate parametrization of an
arc of a curve, via a real rational cubic interpolating 8 points on
the arc.
It remained to find conditions that guaranteed the existence of a
solution which is close to the arc.
This was just solved by Fiedler-Le Touz\'e, an MSRI postdoc who had
studied cubics (not necessarily rational) interpolating 8 points to
help classify real plane curves of degree 9.


\end{document}